\documentclass[12pt]{article}%

\usepackage{amsfonts, amsmath, amscd}
\usepackage[psamsfonts]{amssymb}

\newcommand{\pr}{\par {\bf Proof.} }   
\newcounter{df}
\newcounter{pro}
\newenvironment{pro}{\par
\refstepcounter{pro}
{\bf Proposition \arabic{pro}.} }{}
\newcounter{rem}
\newenvironment{rem}{\par
\refstepcounter{rem}
{\bf Remark \arabic{rem}.} }{}
\newcounter{exa}
\newenvironment{exa}{\par
\refstepcounter{exa}
{\bf Example \arabic{exa}.} }{}
\newcounter{teo}
\newenvironment{teo}{\par
\refstepcounter{teo}
{\bf Theorem \arabic{teo}.} \it }{}
\newcounter{cor}
\newenvironment{cor}{\par
\refstepcounter{cor}
{\bf Corollary \arabic{cor}.} \it }{}
\newcounter{st}
\newcounter{lem}
\newenvironment{lem}{\par
\refstepcounter{lem}
{\bf Lemma \arabic{lem}.} \it }{}
\makeatletter                                                   %
\renewcommand{\section}{\@startsection{section}{1}
{\parindent}{3.5ex plus 1ex minus 0.2ex}{2.3ex plus 0.2ex}{\bf}}
\makeatother                                                    %

\begin{document}%

\author{S.S. Gabriyelyan\footnote{The author was partially supported
 by Israel Ministry of Immigrant Absorption}}
\title{On $T$-sequences and characterized subgroups}
\date{}

\makeatletter
\renewcommand{\@makefnmark}{}
\renewcommand{\@makefntext}[1]{\parindent=1em #1}
\makeatother

\maketitle

\vspace{-\baselineskip}

\begin{abstract}
Let $X$ be a compact metrizable abelian group and  $\mathbf{u}=\{ u_n\}$ be a sequence in its dual $X^{\wedge}$. Set $s_{\mathbf{u}} (X)= \{ x: (u_n,x)\to 1\}$ and $\mathbb{T}_0^H = \{ (z_n)\in \mathbb{T}^{\infty} : z_n\to 1 \}$. Let $G$ be a subgroup of $X$. We prove that $G=s_{\mathbf{u}} (X)$ for some $\mathbf{u}$ iff it can be represented as some dually closed subgroup $G_{\mathbf{u}}$ of ${\rm Cl}_X G \times \mathbb{T}_0^H$. In particular, $s_{\mathbf{u}} (X)$ is polishable. Let $\mathbf{u}=\{ u_n\}$ be a $T$-sequence. Denote by $(\widehat{X}, \mathbf{u})$ the group $X^{\wedge}$ equipped  with the finest group topology in which $u_n \to 0$. It is proved that $(\widehat{X}, \mathbf{u})^{\wedge} =G_{\mathbf{u}}$ and  $\mathbf{n} (\widehat{X}, \mathbf{u}) = s_{\mathbf{u}} (X)^{\perp}$. We also prove that the group generated by a Kronecker set can not be characterized.
\end{abstract}

We shall write our abelian groups additively. For a topological group $X$, $\widehat{X}$ denotes the group of all continuous characters on $X$. We denote its dual group by $X^{\wedge}$, i.e. the group $\widehat{X}$ endowed with the compact-open topology. A group $X$ equipped with discrete topology is denoted by $X_d$. Denote by $\mathbf{n}(X) = \cap_{\chi\in \widehat{X}} {\rm ker} \chi$ the von Neumann radical of $X$. $X$ is named {\it Pontryagin reflexive} or {\it reflexive} if the canonical homomorphism $\alpha_X : X\to X^{\wedge\wedge} , x\mapsto (\chi\mapsto (\chi, x))$ is a topological isomorphism. If $H$ is a subgroup of $X$, we denote by $H^{\perp}$ its annihilator. Let $X$ and $Y$ be topological groups and $\varphi : X\to Y$ a continuous homomorphism. We denote by $\varphi^{\ast} : Y^{\wedge} \to X^{\wedge}, \chi\mapsto \chi\circ\varphi$, the dual homomorphism of $\varphi$. Let $A$ be a subset of a group $X$. Set $(2)A=A+A, (n+1)A =(n)A+A$. $\langle A\rangle$ denotes the subgroup generated by $A$.

Let $G$ be a Borel subgroup of a Polish group $X$. $G$ is called polishable if there exists a Polish group topology $\tau$ on $G$ such that the identity map $i : (G, \tau ) \to X, i(g)=g,$ is continuous.
We remaind that a Polish group topology on a Borel subgroup of a Polish group is unique.

Let $X$ be a compact metrizable group and $\mathbf{u} =\{ u_n \}$ a sequence of elements of $\widehat{X}$. We denote by $s_{\mathbf{u}} (X)$ the set of all $x\in X$ such that $(u_n , x)\to 1$. Let $G$ be a subgroup of $X$. If $G=s_{\mathbf{u}} (X)$ we say that $\mathbf{u}$ {\it characterizes} $G$ and that $G$ is {\it characterized} (by $\mathbf{u}$) or {\it basic $\mathfrak{g}$-closed}. Denote by $\mathfrak{g}_X (G) := \bigcap \{ s_{\mathbf{u}} (X) :\; \mathbf{u}\in \widehat{X}^{\mathbb{N}}, \; G\leqslant s_{\mathbf{u}} (X) \}$. $G$ is called $\mathfrak{g}$-{\it closed} (resp. $\mathfrak{g}$-{\it dense}) if $G=\mathfrak{g}_X (G)$ (resp. $\mathfrak{g}_X (G) =X$) \cite{DMT}.

The following group plays the key role in our considerations. Set \cite{Gab}
$$\mathbb{T}^H_0 := \left\{ \omega =(z_n )\in \mathbb{T}^{\infty} :  \; z_n \to 1 \right\}, d_0 (\omega_1 , \omega_2 )= \sup \{ |z_n^1 -z_n^2 |, n\in \mathbb{N} \}.
$$
Then $\mathbb{T}^H_0$ is a Polish group and $\left( \mathbb{T}^H_0\right)^{\wedge} = \mathbb{Z}^{\infty}_0 =\{ \textbf{n}= (n_1,\dots, n_k, 0,\dots) , n_j \in \mathbb{Z} \}.$ The topology on $\mathbb{Z}^{\infty}_0$ is described in \cite{Gab}. Note that $\mathbb{T}^H_0$ is reflexive and characterized subgroup of $\mathbb{T}^{\infty}$ by the sequence $e_1 =(1,0,0,\dots), e_2 = (0,1,0,\dots), \dots$

The article is divided on two parts.

{\bf I. Characterized subgroups.} Let $X=\mathbb{T}$ and $\mathbf{u} =\{ u_n \}$ be a sequence of integers. Then the elements of $s_{\mathbf{u}} (\mathbb{T})$ are called the {\it topologically $\mathbf{u}$-torsion elements} of $\mathbb{T}$. The interest in topologically $\mathbf{u}$-torsion elements is determined, in particular, by the following three arguments. The first, this notion generalized the algebraic notion of torsion element, naturally connecting the algebraic and the topological structures. The second, it was shown in \cite{BDM}, \cite{Rac} that the subgroups of the form $s_{\mathbf{u}} (\mathbb{T})$ lead to the description of precompact group topologies with converging sequences (see also \cite{BDW}, \cite{DMT}). The third, there are some applications to Diophantine approximation, dynamical systems and ergodic theory \cite{PeS}, \cite{Win}. The study of subgroups $s_{\mathbf{u}} (\mathbb{T})$ is interesting in two directions: {\it concrete} and {\it general}.  The concrete one is the description of $s_{\mathbf{u}} (\mathbb{T})$ either for some {\it concrete} $\mathbf{u}$ \cite{Arm}, \cite{BD2}, \cite{DiK}, or for sequences $\mathbf{u}$ that verify a {\it linear recurrence relation} \cite{BD2}, \cite{DDS}, \cite{KrL}, \cite{Lar}. The general directions means the finding  a general criterion for deciding whether a subgroup $H$ is characterizable. The answer to this question is not simple even for a cyclic countable subgroup of $\mathbb{T}$. It was proved in \cite{Lar} that any cyclic subgroup of $\mathbb{T}$, generated by an irrational number with bounded continued fraction coefficients, is characterizable. This was extended to all countable subgroups of the circle group in \cite{BDS}.
In the general case when $X$ is an infinite metrizable compact group and $\mathbf{u}$ is a sequence in $X^{\wedge}$, it is known very little about characterized subgroups. D.~Dikranjan and K.~Kunen \cite{DiK} proved that every countable infinite subgroup is characterizable and found a condition under which the countable union of an increasing family of compact groups is characterizable. In the following theorem we describe {\it all} characterized subgroups as follows.

\begin{teo} \label{t1}
{\it Let $G$ be a subgroup of an infinite compact metrizable abelian group $X$.
\begin{enumerate}
\item Let $G$ be dense in $X$. Then $G$ is characterized by $\mathbf{u}$ if and only if there exists a dually closed subgroup $G_{\mathbf{u}}$ of $X\times \mathbb{T}^H_0$ such that the restriction to $G_{\mathbf{u}}$ of the projection $\pi_X , \pi_X (x,\omega) =x,$ is a bijection onto $G$ and for every $n$ the character $u_n \in X^{\wedge}$ satisfies
$$
(u_n , g ) = \pi_n (g,\omega ),\; \forall (g,\omega )\in G_{\mathbf{u}},
$$
where $\pi_n (x, (z_n)) = z_n$.
\item $G$ is characterized in $X$ if and only if it is characterized in ${\rm Cl}_X (G)$.
\end{enumerate} }
\end{teo}

\begin{cor} \label{c1}
{\it Each characterized subgroup $G = s_{\mathbf{u}} (X)$ admits a locally quasi-convex Polish group topology.}
\end{cor}

We denote $G = s_{\mathbf{u}} (X)$  with this topology by $G_{\mathbf{u}}$ too.

Since any countable subgroup is $\sigma$-compact, D.~Dikranjan and K.~Kunen also posed the following question: is any $\sigma$-compact subgroup is characterizable. In general, a negative answer was given by A.~Bir\'{o} \cite{Bir} as a corollary of a non-trivial result of J.~Aaronson and N.~Nadkarni \cite{AaN}. He proved that the group $\langle K \rangle$ generated by an {\it uncountable} Kronecker subset $K$ of $\mathbb{T}$ is not characterizable. Let us remind that
a subset $K$ of an infinite compact metrizable abelian group $X$ is called a Kronecker set if it is nonempty, compact, and for every continuous function $f: K\to \mathbb{T}$ and $\varepsilon >0$ there exists a character $\chi \in X^{\wedge}$ such that
$$
\max \{ \| f(x) - (\chi , x) \|, \; x\in K \} <\varepsilon.
$$
There are many uncountable Kronecker sets in any compact abelian group, which contains no compact open torsion subgroup (theorem 41.5 \cite{HR2}). Since the role of Kronecker sets in constructing curios examples in harmonic analysis is very important (see, for example, $\S$ 41 \cite{HR2} and \cite{GrM}), we may expect that Kronecker sets are pathological from the point of view of topological algebra and this is indeed so. We remark that the pathology of Kronecker sets is defined by its algebraical independency. A subset $E$ of $X$ is called (algebraically) independent if whenever $x_1, \dots, x_m \in E$ and $n_1, \dots, n_m \in  \mathbb{Z}$, and $\sum_{i=1}^m n_i x_i =0$, then $n_i x_i =0$ for $1\leq i\leq m$.   Every  Kronecker set is independent and all its elements have only infinite order (theorem 41.8 \cite{HR2}).
We prove the following.
\begin{teo} \label{t7}
{\it Let $X$ be an infinite compact metrizable group and $K$ be an uncountable Kronecker set. Then $\langle K \rangle$ is not polishable. In particular,  $\langle K \rangle$ can not be characterized.}
\end{teo}

Another example of a non characterized subgroup is defined in \cite{Gab}: the group $\mathbb{T}^H_2$ is a {\it polishable} $\sigma$-compact subgroup of $\mathbb{T}^{\infty}$ but its Polish group topology is {\it not} locally quasi-convex. The next example shows that there exists a $\sigma$-compact non characterized subgroup which admits a {\it reflexive} (hence, {\it locally quasi-convex}) Polish group topology. In other words, corollary \ref{c1} cannot be inverted.

\begin{exa} \label{exa1}
{\it Set {\rm \cite{Gab}}
$$
\mathbb{T}^H_1 := \left\{ \omega =(z_n )\in \mathbb{T}^{\infty} :
\; \sum_{n=1}^{\infty} |1-z_n | <\infty \right\} ,
d(\omega_1 , \omega_2 )= \sum_{n=1}^{\infty} |z_n^1
-z_n^2 |.
$$
Then $\mathbb{T}^H_1$ has the following properties.
\begin{enumerate}
\item It is $\sigma$-compact polishable subgroup of $\mathbb{T}^{\infty}$. Its Polish group topology, which is defined by  the metric $d(\omega_1 , \omega_2 )$, is reflexive.
\item It is $\mathfrak{g}$-closed but not a basic $\mathfrak{g}$-closed subgroup of $\mathbb{T}^{\infty}$.
\end{enumerate} }
\end{exa}

It is known that there exists a {\it Borel} $\mathfrak{g}$-dense subgroups \cite{HaK}. We do not know the answer of the following

\begin{enumerate}
\item[] {\bf Question 1}. Do there exist a proper $\mathfrak{g}$-dense {\it polishable} subgroup of a compact group?
\end{enumerate}

{\bf II. Duality and reflexivity.}
Following E.G.Zelenyuk and I.V.Protasov \cite{ZP1}, \cite{ZP2}, we say that a sequence $\mathbf{u} =\{ u_n \}$ in a group $G$ is a $T$-{\it sequence} if there is a Hausdorff group topology on $G$ for which $u_n $ converges to zero. The group $G$ equipped with the finest group topology with this property is denoted by $(G, \mathbf{u})$.

As it was mentioned above, there is a connection between the characterized subgroups of $X$ and the {\it precompact} group topologies with converging sequences on $\widehat{X}$. But, we show that the characterized subgroups of $X$ are related, as well, to {\it Hausdorff} group topologies with converging sequences on $\widehat{X}$: we prove that if $\mathbf{u}$ is a $T$-sequence in $X^{\wedge}$ then $G_{\mathbf{u}}$ is the dual group of $(\widehat{X}, \mathbf{u})$. Moreover, it is proved that the von Neumann radical $\mathbf{n} (\widehat{X}, \mathbf{u})$ is  $s_{\mathbf{u}} (X)^{\perp}$.

\begin{teo} \label{t2}
{\it Let $G=s_{\mathbf{u}} (X)$ be a dense characterized subgroup of an infinite compact metrizable abelian group $X$. Then:
\begin{enumerate}
\item[1)] $\mathbf{u}$ is a $T$-sequence;
\item[2)] $(\widehat{X}, \mathbf{u}) = \left( X\times \mathbb{T}^H_0 \right)^{\wedge} / G_{\mathbf{u}}^{\perp}$ and the projection $\left( X\times \mathbb{T}^H_0 \right)^{\wedge} \to \left( X\times \mathbb{T}^H_0 \right)^{\wedge} / G_{\mathbf{u}}^{\perp}$ is compact covering;
\item[3)] $(\widehat{X}, \mathbf{u})$ is a complete hemicompact $k$-space;
\item[4)] $(\widehat{X}, \mathbf{u})^{\wedge} = G_{\mathbf{u}}$.
\end{enumerate} }
\end{teo}

\begin{teo} \label{t3}
{\it Let $X$ be an infinite compact metrizable abelian group and $\mathbf{u} =\{ u_n \}$ be a $T$-sequence in $X^{\wedge}$. Then $(\widehat{X}, \mathbf{u})^{\wedge} = G_{\mathbf{u}}$ and $\mathbf{n} (\widehat{X}, \mathbf{u})$ is equal to  $s_{\mathbf{u}} (X)^{\perp}$ algebraically.}
\end{teo}

From this duality we can obtain the following simple characterization of maximal/minimal almost-periodicity of $(H, \mathbf{u})$ for any $T$-sequence in a group $H$.

\begin{cor} \label{c2}
{\it Let $\mathbf{u} =\{ u_n \}$ be a $T$-sequence in a group $H$. Then
\begin{enumerate}
\item $(H, \mathbf{u})$ is maximally almost-periodic if and only if $s_{\mathbf{u}} (H_d^{\wedge})$ is dense in $H_d^{\wedge}$.
\item $(H, \mathbf{u})$ is minimally almost-periodic if and only if $H=\langle \mathbf{u}\rangle$ and $s_{\mathbf{u}} (H_d^{\wedge} )=\{ 0\}$.
\end{enumerate} }
\end{cor}
In particular, item 2 of corollary \ref{c2} yields that a minimally almost-periodic group of the form $(H, \mathbf{u})$ is necessarily countable.

Let us note that $(\mathbb{Z}_0^{\infty} , \mathbf{e}),$ where $\mathbf{e} =\{ e_n\}$, is {\it reflexive} \cite{Gab} (in fact it is a Graev free topological abelian group over the convergent sequence $\mathbf{e}$ and $(\mathbb{Z}_0^{\infty} , \mathbf{e})^{\wedge} = \mathbb{T}^H_0 $). On the other hand, we have the following.
\begin{exa} \label{exa2}
{\it let $X$ be a compact metrizable abelian group and $\mathbf{u}$ be a $T$-sequence in $X^{\wedge}$ such that $s_{\mathbf{u}} (X)$ is countable and dense (it is possible by theorem {\rm 1.4 \cite{DiK}}). Then the bidual group $(\widehat{X}, \mathbf{u})^{\wedge\wedge}$ is compact and $\alpha_{(\widehat{X}, \mathbf{u})} (\widehat{X})$ is its dense subset. In particular, $(\widehat{X}, \mathbf{u})$ is not reflexive.}
\end{exa}

\pr By theorem \ref{t3}, $(\widehat{X}, \mathbf{u})^{\wedge}$ is $s_{\mathbf{u}} (X)$ with its (unique) Polish group topology, that must necessarily be the discrete one. Hence $(\widehat{X}, \mathbf{u})^{\wedge}$ is discrete, and, consequently, $(\widehat{X}, \mathbf{u})^{\wedge\wedge}$ is compact. Since the topology of $(\widehat{X}, \mathbf{u})$ is never bounded (theorem 2.3.12 \cite{ZP2}), $(\widehat{X}, \mathbf{u})$ is not reflexive.

Since $G_{\mathbf{u}}= (\widehat{X}, \mathbf{u})^{\wedge}$ is discrete and dense in a compact group $X$, then, by (24.41) \cite{HR1}, $\widehat{X}$  is dense in the compact group $G_{\mathbf{u}}^{\wedge}$. Hence $\alpha_{(\widehat{X}, \mathbf{u})} (\widehat{X})$ is dense. $\Box$

Hence the following questions are meaningful: assume that $\mathbf{u}$ is a $T$-sequence in an infinite group $H$.

\begin{enumerate}
\item[] {\bf Question 2}. When $(H, \mathbf{u})$ is reflexive?
\item[] {\bf Question 3}. When $(H, \mathbf{u})^{\wedge}$ is reflexive?
\end{enumerate}
Note that the subgroup $\langle\mathbf{u}\rangle$ is open in $(H, \mathbf{u})$. Then, by theorem 2.3 \cite{BCM}, $(H, \mathbf{u})$ is reflexive if and only if $\langle\mathbf{u}\rangle$ is reflexive. By lemma 2.2 and theorem 2.6 \cite{BCM}, $(H, \mathbf{u})^{\wedge}$ is reflexive if and only if $\langle\mathbf{u}\rangle^{\wedge}$ is reflexive. Thus we can restrict our attention to the case when $H$ is countable. Let $H$ be countable and infinite. Taking into account theorem \ref{t2} and \ref{t3}, it is more convenient to assume that $H= \widehat{X}$, where $X=H_d^{\wedge}$ is an infinite compact metrizable group. Thus question 3 we can restate as follows.

\begin{enumerate}
\item[] {\bf Question 3$'$}. When the Polish group $G_{\mathbf{u}}$ is reflexive?
\end{enumerate}

\begin{pro} \label{p1}
{\it Let $\alpha_{(\widehat{X}, \mathbf{u})} (\widehat{X})$ be dense in $G_{\mathbf{u}}^{\wedge}$. Then
\begin{enumerate}
\item[1)]
$(\widehat{X}, \mathbf{u})$ is locally quasi-convex if and only if $(\widehat{X}, \mathbf{u})$ is reflexive. In particular,  $s_{\mathbf{u}} (X)$ is dense in $X$ and $\mathbf{n} (\widehat{X}, \mathbf{u}) =\{ 0\}$.
\item[2)] $G_{\mathbf{u}}$ is reflexive.
\end{enumerate} }
\end{pro}

\begin{teo} \label{t6}
{\it $(\widehat{X}, \mathbf{u})$ is reflexive if and only if the following three conditions satisfy:
\begin{enumerate}
\item $s_{\mathbf{u}} (X)$ is dense in $X$ (and, hence, $\mathbf{u} $ is a $TB$-sequence).
\item $G_{\mathbf{u}}$ is dually embedded in $X\times \mathbb{T}^H_0$.
\item $(\widehat{X}, \mathbf{u})$ is locally quasi-convex.
\end{enumerate} }
\end{teo}

For example, if $s_{\mathbf{u}} (X)$ is countable and dense, then, by example \ref{exa2} and proposition \ref{p1},  $(\widehat{X}, \mathbf{u})$ is {\it not locally quasi-convex}. We do not know the answers on the following questions. Let $X$ be a compact metrizable group and $\mathbf{u}$ a $T$ sequence in $X^{\wedge}$.

\begin{enumerate}
\item[] {\bf Question 4}. Under which conditions on $\mathbf{u}$ the group $(\widehat{X}, \mathbf{u})$ is (non) locally quasi-convex?
\item[] {\bf Question 5}. Is there a sequence $\mathbf{u}$ such that $G_{\mathbf{u}}$ is not reflexive?
\end{enumerate}
Note that if $G_{\mathbf{u}}$ is dually embedded in $X\times \mathbb{T}^H_0$, then it is reflexive by lemma 14.6 \cite{Ban} and corollary 3 \cite{Cha}.

Our interest to questions 3 and 5 is motivated also by the following. If any $G_{\mathbf{u}}$ is reflexive, we obtain a natural class of reflexive groups which is, on the one hand, more wide than the class of locally compact groups. On the other hand, this class is not ``very big'' and can be transferred to the non abelian case (see discussion of this problem in \cite{Gab}).

Assume that there exists a $T$-sequence $\mathbf{u}$ such that $G_{\mathbf{u}}$ is not reflexive. Since $(\widehat{X}, \mathbf{u})$ is a $k$-space and $G_{\mathbf{u}}$ is metrizable, then $\alpha_{(\widehat{X}, \mathbf{u})}$ and $\alpha_{G_{\mathbf{u}}}$ are continuous (corollary 5.12 \cite{Aus}). Thus, by the Vilenkin theorem \cite{Vil}, \cite{Aus}, we must have $G_{\mathbf{u}}^{\wedge\wedge} = G_{\mathbf{u}}\oplus Y$ and $(\widehat{X}, \mathbf{u})^{\wedge\wedge\wedge} = (\widehat{X}, \mathbf{u})^{\wedge} \oplus Y$, for some Polish group $Y$. We are not aware of any group with these properties.

\begin{center}
{\bf\large The proofs}
\end{center}

{\bf Proof of theorem \ref{t1}}. 1. Assume that $G$ is dense.

$\Rightarrow$ Let $G= s_{\mathbf{u}} (X)$. Set $G_{\mathbf{u}} =\left\{ \left(g; (u_1 , g),(u_2 , g),\dots \right), \; g\in G\right\}$.

Let us prove that $G_{\mathbf{u}}$ is dually closed (and, hence, is closed). Let $\mathbf{x} = (x, (z_n) ) \not\in G_{\mathbf{u}}$. Then there exists the minimal index $i$ such that $z_i \not= (u_i , x)$. Let
$$
p_i : X\times \mathbb{T}^H_0 \to X\times \mathbb{T}^i , p_i (x, (z_n))= (x, z_1, \dots , z_i),
$$
be the natural projection and $L_i = \{ (x, (u_1 , x), \dots , (u_i , x)), x \in X\}$. Then $L_i$ is a closed subgroup of $X\times \mathbb{T}^i$, $p_i (G_{\mathbf{u}})\subset L_i$ and $p_i (\mathbf{x})\not\in L_i$. Since $X\times \mathbb{T}^i$ is compact, there exists its character $\omega_i = (u, n_1,\dots , n_i)$ such that $(\omega_i , L_i )=1$, but $(\omega_i , p_i (\mathbf{x}))\not= 1$. Set $\widetilde{\omega} =(u, n_1,\dots , n_i , 0 ,\dots) \in (X\times \mathbb{T}^H_0 )^{\wedge}$. It is clear that
$$
(\widetilde{\omega} ,G_{\mathbf{u}}) \subset (\omega_i ,L_i ) = \{ 1\} , \mbox{ but } (\widetilde{\omega} , \mathbf{x} ) = (\omega_i , p_i (\mathbf{x}))\not= 1.
$$
Hence $G_{\mathbf{u}}$ is dually closed.

Since $X\times \mathbb{T}^H_0$ is reflexive, its topology is locally quasi-convex. Thus the induced topology on $G_{\mathbf{u}}$ is locally quasi-convex too. By the construction, $(u_n , g)=\pi_n (g,\omega)$ and the restriction of $\pi_X$ to $G_{\mathbf{u}}$ is bijective.

$\Leftarrow$ Denote $\pi_X (G_{\mathbf{u}})$ by $G$ and set $\mathbf{u} =\{ u_n \}$. Let us prove that $G=s_{\bf u} (X)$. It is clear that $(u_n , g)\to 1 , \forall g\in G$. Hence $G\subset s_{\bf u} (X)$. Let us prove the converse inclusion. If $x\in s_{\bf u} (X)\setminus G$, then $\mathbf{x} := (x, (u_n, x)) \not\in G_{\mathbf{u}}$. Since $G_{\mathbf{u}}$ is dually closed,  there exists $\widetilde{\omega} = (u, l_1, \dots, l_m, 0,\dots)\in (X\times \mathbb{T}^H_0 )^{\wedge}$ such that $\widetilde{\omega}\in G_{\mathbf{u}}^{\perp}$ and $(\widetilde{\omega} , \mathbf{x}) \not= 1$. Therefore
$$
(u,g)\cdot (u_1 , g)^{l_1}\dots (u_m , g)^{l_m} = (u+l_1 u_1 +\dots + l_m u_m , g) =1 , \forall g\in G.
$$
Since $G$ is dense, the last equality is possible only if  $u+l_1 u_1 +\dots + l_m u_m =0$. Hence $(\widetilde{\omega} ,\mathbf{x}) =1$. This is a contradiction.

2. Let $G=s_{\bf u} (X)$. By lemma 2.8 \cite{DiK}, $s_{\bf u} (X)= s_{\bf \widetilde{u}} (Y)$, where $Y= {\rm Cl}_X \left( s_{\bf u} (X)\right)$ and $\widetilde{u}_n = u_n|_Y$. Hence $G$ is characterized in ${\rm Cl}_X (G)$.

The inverse assertion follows from lemma 2.1 \cite{DiK}.  $\Box$

{\bf Proof of corollary \ref{c1}}. By lemma 2.8 \cite{DiK}, $s_{\bf u} (X) = s_{\bf \widetilde{u}} (Y)$ where $Y= {\rm Cl}_X (s_{\bf u} (X))$ and $\widetilde{u}_n = u_n |_Y$. By theorem \ref{t1}, $s_{\bf u} (X)$ is polishable. By the reflexivity of the group $X\times \mathbb{T}_0^H$ \cite{Gab}, that group is locally quasi-convex. Thus, as a closed subgroup of a locally quasi-convex group, also $G_{\bf u}$ is is locally quasi-convex by proposition 6.8 \cite{Aus}. $\Box$

{\bf Proof of theorem \ref{t7}}. Let us assume the converse and $r$ be a complete metric on $G=\cup (n)(K\cup (-K))$. The $\varepsilon$-neighborhood of the neutral element we denote by  $U_{\varepsilon}$. By Baire's category theorem, there exists $n_0$ such that $(n_0)(K\cup (-K))$ contains some $U_{\varepsilon}$. Since $G$ is uncountable,  $U_{\frac{\varepsilon}{2n_0}}$ contains a non-zero element $g$. Hence $2n_0 g\in U_{\varepsilon}$. On the other hand, we can represent $g$ and $2n_0 g$ in the form
$$
g= a_{11} g_{11} +\dots + a_{{l_1}1} g_{{l_1}1} , \quad 2n_0 g =a_{12} g_{12} +\dots + a_{{l_2}2} g_{{l_2}2},
$$
$$
\mbox{where } g_{ij} \in K \mbox{ and } |a_{11}| +\dots +|a_{{l_1}1}| \leq n_0, \; |a_{12}| +\dots +|a_{{l_2}2}| \leq n_0.
$$
We can assume that $g_{i1}=g_{i2}, i=1,\dots,s$, and $g_{ij} \not= g_{kq}$ for the remaining combinations of $i,j,k,q$. Since $2n_0 a_{11} g_{11} +\dots +2n_0 a_{{l_1}1} g_{{l_1}1} =a_{12} g_{12} +\dots + a_{{l_2}2} g_{{l_2}2}$, then
$$
\sum_{i=1}^s (2n_0 a_{i1} - a_{i2})g_{i1} +\sum_{i=1}^{l_1 -s} 2n_0 a_{s+i,1} g_{s+i,1} - \sum_{i=1}^{l_2 -s} a_{s+i,2} g_{s+i,2} =0.
$$
Since $K$ is algebraically independent and contains only elements of infinite order (theorem 41.8 \cite{HR2}),  $a_{s+i,1} =a_{s+j,2}=2n_0 a_{i1} - a_{i2} =0$. This possible only if $a_{i1}= a_{j2}=0 ,\forall i,j$. Hence $g=0$. It is a contradiction.  $\Box$

{\bf Proof of example \ref{exa1}}. The items 1 was proved in \cite{Gab}. Let us prove item 2.

Let $\mathbf{\omega} = \{ \omega_k \}$ be a $T$-sequence in $(\mathbb{T}^{\infty})^{\wedge} = \mathbb{Z}^{\infty}_0$ such that $\mathbb{T}^H_1 \subset s_{\mathbf{\omega}} (\mathbb{T}^{\infty}) =s_{\mathbf{\omega}}$. Suppose that
$$
\omega_k = (0,\dots,0, n^k_{r_{1}^k}, \dots, n^k_{r_{s_k}^k}, 0,\dots) \in \mathbb{Z}_0^{\infty}, \; n^k_{r_{i}^k} \not= 0 \mbox{ for } 1\leq i \leq s_k,
$$
where $n^k_{r_{i}^k}$ occupies the position $r_{i}^k$.

1) {\it Let us prove that $r_1^k \to\infty$.}

Assume the converse and there exists an index $a$ such that $r_{1}^{k_j} =a$ for $k_1 < k_2 <\dots$ Denote by $\pi_a$ and $\pi'_a$ the projections of $\mathbb{T}^{\infty}$ and  $\mathbb{Z}_0^{\infty}$ on the coordinate $a$. Hence the sequence $\mathbf{u} = \{ \pi'_a (\omega_k) \}$ is non trivial. Since any $\mathbf{z} =(z_n)$, where $z_a =z$ and $z_n =1$ if $n\not= a$, belongs to $\mathbb{T}^H_1$, we have $s_{\mathbf{u}} (\mathbb{T}) =\mathbb{T}$. By corollary \ref{c3} (see below),
$\mathbf{u}$ is trivial. It is a contradiction.

2) {\it Set $d_k =\max \{ |n^k_{r_{i}^k} |, i=1,\dots, s_k \}$. Let us prove that there exists $C>0$ such that $d_k <C$ for every $k$.}

Assume the converse and $d_{k_j} \to\infty$ at $j\to\infty$. By 1) we can suppose that
\begin{equation} \label{e1}
d_{k_j} = |n^{k_j}_{r_{i_j}^{k_j}} | > j^2 \; \mbox{ and } \; r_{s_{k_j}}^{k_j} < r_{1}^{k_{j+1}}.
\end{equation}

Set $\mathbf{z} = (z_n)$, where
$$
z_l =\exp \left\{ 2\pi i \cdot \frac{1}{2d_{k_j}} \right\} , \mbox{ if } l= r_{i_j}^{k_j} , \mbox{ and } z_l =1 \mbox{ otherwise}.
$$
Then, by (\ref{e1}),
$$
\sum_l | 1-z_l | = \sum_j | 1 - \exp \left\{ 2\pi i \cdot \frac{1}{2d_{k_j}} \right\} | < 2\sum_j \frac{1}{j^2} <\infty
$$
and  $\mathbf{z}\in \mathbb{T}^H_1$. But, by (\ref{e1}), we have
$$
(\omega_{k_j} , \mathbf{z}) = \exp \left\{ 2\pi i \cdot \frac{1}{2d_{k_j}} \cdot d_{k_j} \right\} =-1 \not\to 1.
$$
It is a contradiction.

3) {\it Let us prove that $\mathbb{T}^H_1 \not= s_{\mathbf{\omega}}$ and, hence, $\mathbb{T}^H_1$ is not characterizable.}

Set $k_1 =1$. By 1), we can choose $k_2 >k_1$ such that
$$
r_{s_{k_1}}^{k_1} < r_1^k \mbox{ for every } k\geq k_2.
$$
By 1), we can choose $k_{m+1} > k_m$ such that
$$
\max \{  r_{s_{1}}^{1}, r_{s_{2}}^{2}, \dots , r_{s_{k_m}}^{k_m} \} < r_1^k \mbox{ for every } k\geq k_{m+1},
$$
and so on. Set $\mathbf{z} = (z_l)$, where
$$
z_l =\exp \left\{ \frac{2\pi i}{m\cdot C} \right\} , \mbox{ if } l= r_{1}^{k_m} , \mbox{ and } z_l =1 \mbox{ otherwise}.
$$
Then, by the definition of the sequence $k_m$, we have the following: if $m\geq 2$ and $k_m \leq t< k_{m+1}$, there exist three possibilities
\begin{enumerate}
\item[a)] $(\omega_t , \mathbf{z}) =1$;
\item[b)] $(\omega_t , \mathbf{z}) = \exp \left\{ \frac{2\pi i}{m\cdot C} \cdot  n_{r_1^{k_m}}^{t} \right\} $;
\item[c)] $(\omega_t , \mathbf{z}) = \exp \left\{ \frac{2\pi i}{m\cdot C} \cdot n_{r_1^{k_m}}^{t} \right\} \cdot \exp \left\{ \frac{2\pi i}{(m+1)\cdot C} \cdot n_{r_1^{k_{m+1}}}^{t} \right\}$.
\end{enumerate}
Hence $(\omega_t , \mathbf{z}) \to 1$ and  $\mathbf{z} \in s_{\mathbf{\omega}}$. On the other hand, the series
$\sum_l |1-z_l | = \sum_m | 1-\exp \left\{ 2\pi i \cdot \frac{1}{mC} \right\} -1 |$ diverges. Thus $\mathbf{z} \not\in \mathbb{T}^H_1$.

4) {\it Let us prove that $\mathbb{T}^H_1$ is $\mathfrak{g}$-closed.}

Denote by $Y$ the $\mathfrak{g}$-closure of $\mathbb{T}^H_1$ in $\mathbb{T}^{\infty}$. Let $\mathbf{z} =(z_n)$ be such that $\sum_n | 1-z_n | = \infty$. It is enough to prove that $\mathbf{z} \not\in Y$.

(a) If $z_n \not\to 1$, then $(e_n, \mathbf{z}) =z_n \not\to 1$, but $s_{\mathbf{e}} (\mathbb{T}^{\infty}) = \mathbb{T}^H_0$. Thus $\mathbf{z} \not\in Y$.

(b) Let $z_n =\exp \{ 2\pi i\cdot \varepsilon_n \} , 0\leq \varepsilon_n <1,$ and $z_n \to 1$. Choose $k_0$ such that $0\leq \varepsilon_n < 0,01$ for every $n>k_0$.

Let $k_1 > k_0$ be the first integer such that $\frac{1}{3} < \sum_{j=k_0}^{k_1} \varepsilon_j <\frac{1}{2}$. Let $k_{m+1} > k_m$ be the first integer such that
\begin{equation} \label{e2}
\frac{1}{3} < \sum_{j=k_m +1}^{k_{m+1}} \varepsilon_j <\frac{1}{2}
\end{equation}
and so on. Set $\omega_m = (n_k^m)$, where $n_k^m =1$ if $k_m <k \leq k_{m+1}$ and $n_k^m =0$ otherwise. Then
$$
(\omega_m , \mathbf{z}) = z_{k_m +1} \cdot z_{k_m +2}\cdot \dots \cdot z_{k_{m +1}} = \exp \left\{ 2\pi i \sum_{j=k_m +1}^{k_{m+1}} \varepsilon_j \right\}.
$$
Hence, by (\ref{e2}), $(\omega_m , \mathbf{z}) \not\to 1$ and $\mathbf{z} \not\in s_{\mathbf{\omega}}$.

On the other hand, if $\mathbf{z} =\{ z_n \} \in \mathbb{T}^H_1$, then
$$
|(\omega_m , \mathbf{z}) -1| = | z_{k_m +1} \cdot z_{k_m +2}\cdot \dots \cdot z_{k_{m +1}} -1| \leq \sum_{j=k_m +1}^{k_{m+1}} | z_j -1| \to 0.
$$
Hence $\mathbb{T}^H_1 \subset s_{\mathbf{\omega}}$. Thus $z\not\in Y$ and $\mathbb{T}^H_1$ is $\mathfrak{g}$-closed. $\Box$

In the sequel we need some notations. For every sequence $\{ M_i \}$ of subsets of $\widehat{X}$ we put
$$
\sum_i M_{i}^{\ast} = \cup_i \left( M_{0}^{\ast} +\dots + M_{i}^{\ast} \right), \mbox{ where } A^{\ast} = A\cup (-A)\cup \{ 0\} .
$$
For a $T$-sequence $\{ u_n \}$ and $l,m\in \mathbb{N}$, one puts \cite{ZP1}: $A_m = \{ u_n : n\geq m\}$ and
$$
A(k,m) =\left\{ n_1 u_{r_1} +\dots +n_s u_{r_s} | m\leq r_1 <\dots < r_s , n_i \in \mathbb{Z} \setminus \{ 0\} , \sum_{i=1}^s | n_i | \leq k \right\}.
$$

We need the following lemma.
\begin{lem} \label{l1}
{\it Let $\mathbf{u}$ be a $T$-sequence in an abelian group $H$. Then $(H, \mathbf{u})^{\wedge}$, as a subgroup of $H_d^{\wedge}$, coincides with $s_{\mathbf{u}} (H_d^{\wedge})$ and, hence, $\mathbf{n} (H,\mathbf{u})$ is equal to $s_{\mathbf{u}} (H_d^{\wedge})^{\perp}$ algebraically.}
\end{lem}

\pr Since $H_d \to (H, \mathbf{u})$ is a continuous isomorphism, its conjugate $(H, \mathbf{u})^{\wedge}\to H_d^{\wedge}$ is a continuous monomorphism. If $x\in (H, \mathbf{u})^{\wedge}$, then, by the definition of $(H, \mathbf{u})$, $(u_n, x)\to 1$. Thus $x\in s_{\mathbf{u}} (H_d^{\wedge})$. Let us prove that each $x\in s_{\mathbf{u}} (H_d^{\wedge})$ is continuous on $(H, \mathbf{u})$.

Let $\varepsilon >0$. Choose a sequence $n_0 < n_1 <\dots$ such that
$$
| (u_n, x) - 1 | < \frac{\varepsilon}{2^{k+1}}, \quad \forall n \geq n_k.
$$
Set $U= \sum_i A_{n_i}^{\ast}$. Then $U$ is open in $(H, \mathbf{u})$ (theorem 2.1.3 \cite{ZP2}). Let $y\in U$. Then there exists $i$ such that $y= u_{l_0}^{\ast} + u_{l_1}^{\ast} + \dots + u_{l_i}^{\ast}$, where $u_{l_i}^{\ast} \in A_{n_i}^{\ast}, l_i \geq n_i$. Hence
$$
| (y,x) -1|= | (u_{l_0}^{\ast}, x)\dots (u_{l_i}^{\ast} , x) -1| \leq \sum_{k=0}^i |(u_{l_k}^{\ast} , x) -1| <\varepsilon.
$$
Thus $x$ is continuous. $\Box$

{\bf Proof of theorem \ref{t2}.} 1) At the beginning, we prove that {\it there exists a closed subgroup $H\subset X^{\wedge}\times \mathbb{Z}_0^{\infty}$ such that $H^{\perp} =G_{\mathbf{u}}$ and the canonical homomorphism
$$
Q: X^{\wedge} \to (X^{\wedge}\times \mathbb{Z}_0^{\infty}) /H , \; id (y) = ( y; 0,0,\dots),
$$
and its conjugate
$$
Q^{\ast} : \left( (X^{\wedge}\times \mathbb{Z}_0^{\infty}) /H \right)^{\wedge} \to H^{\perp} =G_{\mathbf{u}}
$$
are continuous isomorphisms}.

Since $G_{\mathbf{u}}$ is dually closed subgroup of the reflexive group $X\times \mathbb{T}_0^H$, we have $G_{\mathbf{u}}=G_{\mathbf{u}}^{\perp\perp}$ and there exists a closed subgroup $H\subset X^{\wedge}\times \mathbb{Z}_0^{\infty}$ such that $H^{\perp} =G_{\mathbf{u}}$ ($H:=G_{\mathbf{u}}^{\perp}$) (14.2, \cite{Ban}). Then $Q^{\ast}$ is a continuous isomorphism.

Let us prove that $Q$ is  a continuous isomorphism. Continuity of $Q$ is obvious. It remains to check that $Q$ is an algebraic isomorphism. By definition, $\mathbf{y} = (y; s_1,\dots,s_l,0,\dots)\in G_{\mathbf{u}}^{\perp}=H$ if and only if for any $g\in G$
$$
1= (\mathbf{y}, (g; (u_k , g)))= (y,g) \prod_{k=1}^l (u_k , g)^{s_k} =(y+\sum_{k=1}^l s_k u_k , g)
$$
Since $G$ is dense, we have $y+\sum_{k=1}^l s_k u_k =0$. Therefore every $\mathbf{y} = (y; s_1,\dots,s_l,0,\dots)$ is equivalent to $\mathbf{y}_0 = (y+\sum_{k=1}^l s_k u_k ; 0,0,\dots)$ modulo $H$. Hence $Q$ is surjective. Evidently, if $y_1 \not= y_2$, then $Q(y_1)\not= Q(y_2)$. Therefore $Q$ is an isomorphism.

2) {\it Let us prove that $u_n \to 0$ in the topology induced from $(X^{\wedge}\times \mathbb{Z}_0^{\infty}) /H$.} Thus, since $(X^{\wedge}\times \mathbb{Z}_0^{\infty}) /H$ is Hausdorff, we will prove that {\it the sequence $\mathbf{u}$ is a $T$-sequence}.

If $\widetilde{V}$ is an open neighborhood of the neutral element in $(X^{\wedge}\times \mathbb{Z}_0^{\infty}) /H$, then there exists an open neighborhood $V$ of zero in $\mathbb{Z}_0^{\infty}$ such that $\pi (\{ 0\} \times V) \subset \widetilde{V}$, where $\pi : X^{\wedge}\times \mathbb{Z}_0^{\infty} \to (X^{\wedge}\times \mathbb{Z}_0^{\infty}) /H$ is the canonical homomorphism. Set
$$
F_{\varepsilon}^l := \left\{ \chi =(0,\dots , 0, n_{l+1}, \dots , n_m , 0,\dots ) : \; n_{l+1}^2 +\dots + n_m^2 \leq \frac{1}{\varepsilon^2} \right\} \subset \mathbb{Z}_0^{\infty}.
$$
Then there exists $n_0$ such that $F_1^{n_0} \subset V$ \cite{Gab}. Since
\begin{equation} \label{1}
(u_n ; 0,0,\dots) \sim (0; (0,\dots, -1, 0, \dots)),
\end{equation}
where $-1$ occupies position $n$, we have $u_n \in \widetilde{V}$ for all $n\geq n_0$. Hence $u_n \to 0$.

3) {\it Let us prove that the topology of $(X^{\wedge}\times \mathbb{Z}_0^{\infty}) /H$ is the finest topology such that $u_n \to 0$.} Thus, by the definition of $T$-sequence, we will prove that $(\widehat{X}, \mathbf{u})= \left( X\times \mathbb{T}_0^H \right)^{\wedge} / G_{\mathbf{u}}^{\perp}$.

Let $\tau''$ be a Hausdorff group topology such that $u_n \to 0$ and $U\in \tau''$ be a symmetric neighborhood of the neutral element and $y\in U$. Choose a sequence $\{ V_k \}$ and $W$ of symmetric neighborhoods of the neutral element such that
\begin{equation} \label{2}
y +W \subset U, \quad U \supset W \supset (2)V_1  \supset (4)V_2 \dots \supset (2^k)V_k \supset\dots
\end{equation}
For any $k$ we can choose $n_k$ such that $u_n \in V_k$ for every $n\geq n_k$. Let $\varepsilon >0$. Choose $k_0$ such that $2^{k_0} > \frac{1}{\varepsilon^3}$. If $\widetilde{\omega} =(0; (0,\dots, 0, m_{n_{k_0} +1}, \dots, m_s, 0, \dots)) \in \{ 0\} \times F_{\varepsilon}^{n_{k_0}}$, then the number of nonzero $m_i$ is at most $\frac{1}{\varepsilon^2}$. Thus
$$
|m_{n_{k_0} +1} | +\dots +|m_s| \leq \frac{1}{\varepsilon^2} \sqrt{m_{n_{k_0} +1}^2 +\dots +m_s^2} <\frac{1}{\varepsilon^3} < 2^{k_0}.
$$
Therefore we have
$$
\left( Q^{-1} \circ\pi\right) (\widetilde{\omega}) = m_{n_{k_0} +1} u_{n_{k_0}+1} +\dots + m_s u_s \in (|m_{n_{k_0} +1} | +\dots +|m_s| ) V_{k_0} \subset (2^{k_0}) V_{k_0} \subset W.
$$
Thus, for any $\mathbf{y}$ such that $\pi(\mathbf{y}) = Q(y)$, we have $\left( Q^{-1} \circ\pi\right) \left( \mathbf{y} + \left( \{ 0\} \times F_{\varepsilon}^{n_{k_0}} \right) \right) \subset y+ W \subset U$. Hence $\mathbf{y} + \left( \{ 0\} \times F_{\varepsilon}^{n_{k_0}} \right) \subset \pi^{-1} ( Q(U))$. As it was proved in item 6 of theorem 1 \cite{Gab},  $\pi^{-1} ( Q(U))$ is open in $X^{\wedge}\times \mathbb{Z}_0^{\infty}$. Since $\pi$ is surjective and open, then $\pi (\pi^{-1} ( Q(U))) = Q(U)$ is open in $(X^{\wedge}\times \mathbb{Z}_0^{\infty}) /H$. Hence the topology on $(X^{\wedge}\times \mathbb{Z}_0^{\infty}) /H$ is finest such that $u_n \to 0$..

4) {\it Let us prove that $\pi$ is compact covering}.

 Let $K$ be a compact subset of $(\widehat{X}, \mathbf{u})$. Since $\langle \mathbf{u} \rangle$ is an open subgroup in $(\widehat{X}, \mathbf{u})$ and $\widehat{X}$ is countable, we can assume that $K\subset \langle \mathbf{u} \rangle$. Let us prove that there exists $k>0$ such that
$$
K\subset A(k,0).
$$
Assume the converse and there exists a sequence $\{ y_n \} \subset K$  such that
$$
y_n \in A(k_n ,0)\setminus A(k_n -1 , 0), \quad k_1 < k_2 <\dots
$$
Since $K$ is compact, there exists a cluster point $y$ of $\{ y_n \}$. Since $K\subset \langle \mathbf{u} \rangle$, there exists $k_0$ such that $y\in A(k_0 ,0)$. Hence we can assume that $y=0\in K$. By lemma 2.3.2 \cite{ZP2}, there exists a neighborhood $U$ of zero such that the set $\{ y_n \} \cap U= \emptyset$. Hence $y$ is not a cluster point. It is a contradiction.

By definition, $A(k,0) \subset \pi (\{ 0 \} \times F_{\frac{1}{k^2}}^0 )$. Since $F_{\varepsilon}^{l}$ is compact in $\mathbb{Z}_0^{\infty}$ \cite{Gab}, $\pi$ is compact covering.

5) {\it Now let us prove that $(\widehat{X}, \mathbf{u})$ is a complete hemicompact $k$-space and $(\widehat{X}, \mathbf{u})^{\wedge} = G_{\mathbf{u}}$}.

Since $X^{\wedge} \times \mathbb{Z}_0^{\infty}$ is a hemicompact $k$-space,   $(\widehat{X}, \mathbf{u})$ is a hemicompact $k$-space by theorem 3.3.23 \cite{Eng}. By theorem 2.3.11 \cite{ZP2}, $(\widehat{X}, \mathbf{u})$ is complete.

By lemma \ref{l1}, $\pi^{\ast} (X^{\wedge} ) = G_{\mathbf{u}}$. Since $\pi : X^{\wedge} \times \mathbb{Z}_0^{\infty} \to (\widehat{X}, \mathbf{u})$ is compact covering, by lemma 5.17 \cite{Aus}, $\pi^{\ast} : (\widehat{X}, \mathbf{u})^{\wedge} \to X\times \mathbb{T}_0^{\infty}$ is an embedding. Hence the corestriction  $(\widehat{X}, \mathbf{u})^{\wedge} \to G_{\mathbf{u}}$ is an open continuous isomorphism. Thus $(\widehat{X}, \mathbf{u})^{\wedge} = G_{\mathbf{u}}$. $\Box$

The following corollary is contained in lemma 3.4.1 \cite{Dik} (see also example 2.8 \cite{BDM}).

\begin{cor} \label{c3}
{\it Let $\mathbf{u} =\{ u_n\}$ be a sequence in $X^{\wedge}$. Then $s_{\mathbf{u}} (X) = X$ if and only if $\mathbf{u}$ is trivial.}
\end{cor}

\pr Let $s_{\mathbf{u}} (X) = X$. By theorem \ref{t2},
$\mathbf{u}$ is a $T$-sequence. Since $(\widehat{X}, \mathbf{u})$ is a $k$-space and a MAP, $\alpha_{(\widehat{X}, \mathbf{u})}$ is a continuous monomorphism. By theorem \ref{t2},
$(\widehat{X}, \mathbf{u})^{\wedge\wedge} = X^{\wedge}$ is discrete, then $(\widehat{X}, \mathbf{u})$ must be discrete too. Hence $\mathbf{u}$ is trivial. The converse is evidently. $\Box$

For the countable case, we can prove theorem  4.20 from \cite{DCl}.

\begin{cor} \label{c4}
{\it  $H\leqslant X$ is $\mathfrak{g}$-dense if and only if $(\widehat{X}, \tau_H )$ has no nontrivial convergent sequences.}
\end{cor}

\pr Let $H$ be $\mathfrak{g}$-dense. Since every closed subgroup of $X$ is characterizable, $H$ is dense in $X$. Let $\mathbf{u} =\{ u_n \}$ be a non trivial convergent sequence in $(\widehat{X}, \tau_H )$. We can assume that $u_n \to 0$ and, hence, $(u_n ,x)\to 1$  for every $x\in H$. By theorem  \ref{t2} and corollary \ref{c3}, $(\widehat{X}, \mathbf{u})^{\wedge} = s_{\mathbf{u}} (X) \not= X.$ Hence $H$ is not $\mathfrak{g}$-dense. It is a contradiction.

Let us prove the converse. Let $\mathbf{u}$ be a sequence in $\widehat{X}$ such that $H\leqslant s_{\mathbf{u}} (X)$. Since $(\widehat{X}, \tau_H )$ is MAP, then $H$ is dense in $X$. By theorem \ref{t2}, $\mathbf{u}$ is a $T$-sequence. By the definition of $\tau_H$ and theorem \ref{t2},
the identity map $(\widehat{X}, \mathbf{u})\to (\widehat{X}, \tau_H )$ is continuous. Hence $u_n \to 0$ in $(\widehat{X}, \tau_H )$. Thus $\mathbf{u}$ is trivial. Hence $s_{\mathbf{u}} (X) = X$. $\Box$

For the proof of theorem \ref{t3}
we need the following lemma which shows that a quotient of a group of the form $(G, \mathbf{u})$ has this form too. We also remark that subgroups of a group $(G, \mathbf{u})$, in general, may have no this form \cite{BaT}.

\begin{lem} \label{l2}
{\it Let $\mathbf{u} =\{ u_n \}$ be a $T$-sequence in a group $G$ and  $H$ be a closed subgroup of $(G, \mathbf{u})$. Let $\pi_0$ be the natural homomorphism from $(G, \mathbf{u})$ to $Y:= (G, \mathbf{u}) /H$. Then $\mathbf{\widetilde{u}} =\{ \widetilde{u}_n \}$, where $\widetilde{u}_n = \pi_0 (u_n)$, is a $T$-sequence in $Y$ and $(Y , \widetilde{\mathbf{u}}) =Y$.}
\end{lem}

\pr Since $H$ is closed, $Y$ is Hausdorff. Since $\widetilde{u}_n =\pi_0 (u_n)\to 0$ in $Y$, $\mathbf{\widetilde{u}}$ is a $T$-sequence.

Let us prove that $(Y ,\widetilde{\mathbf{u}}) =Y$. By definition, the topology $\widetilde{\tau}$ on $(Y ,\widetilde{\mathbf{u}})$ is finer than the topology $\tau$ on $Y$. Let us see now that $\tau$ is finer than $\widetilde{\tau}$. If not, denote by $\mathcal{U}$ and $\mathcal{\widetilde{U}}$ basises at $0$ of $(G, \mathbf{u})$ and $(Y ,\widetilde{\mathbf{u}})$ respectively. Then the sets $U\cap \pi_0^{-1} (\widetilde{U}), U\in \mathcal{U}, \widetilde{U} \in \mathcal{\widetilde{U}}$, form a new Hausdorff topology on $G$ which is finer than the topology of $(G ,\mathbf{u})$ and in which $u_n \to 0$. This is a contradiction. $\Box$

{\bf Proof of theorem  \ref{t3}.}
By lemma \ref{l1}, $\mathbf{n} (\widehat{X}, \mathbf{u}) = s_{\mathbf{u}} (X)^{\perp}$ algebraically.

By definition, $\mathbf{n} (\widehat{X}, \mathbf{u})$ is closed. Let $\pi_0$ be the natural homomorphism from $(\widehat{X}, \mathbf{u})$ to $Y:= (\widehat{X}, \mathbf{u})/ \mathbf{n} (\widehat{X}, \mathbf{u})$. Then, by the definition of $\mathbf{n} (\widehat{X}, \mathbf{u})$, $\pi_0^{\ast} : Y^{\wedge} \to (\widehat{X}, \mathbf{u})^{\wedge}$ is a continuous isomorphism.

Set $\mathbf{\widetilde{u}} =\{ \widetilde{u}_n \}$, where $\widetilde{u}_n = \pi_0 (u_n)$. Then, by lemma \ref{l2}, $(Y ,\widetilde{\mathbf{u}}) =Y$.
Since $\mathbf{n} (\widehat{X}, \mathbf{u}) = s_{\mathbf{u}} (X)^{\perp}$, the group of characters of $Y^{\wedge}_d$ is $H= {\rm Cl}_X (s_{\mathbf{u}} (X))$. Since $ (\widetilde{u}_n , x)= (u_n , x), x\in H$, we have $s_{\mathbf{\widetilde{u}}} (H) =s_{\mathbf{u}} (X)$ is dense in $H$. By theorem \ref{t2},
$Y^{\wedge} = G_{\mathbf{u}}$ is Polish.

On the other hand, by corollary 4.1.5 \cite{ZP2}, $(\widehat{X}, \mathbf{u})$ is a hemicompact $k$-space. Hence $(\widehat{X}, \mathbf{u})^{\wedge}$ is a completely metrizable group \cite{Aus}. Since $Y^{\wedge}$ is separable, then $(\widehat{X}, \mathbf{u})^{\wedge}$ is separable and, hence, Polish. Thus $\pi_0^{\ast}$ is a topological isomorphism and $(\widehat{X}, \mathbf{u})^{\wedge} =Y^{\wedge} =G_{\mathbf{u}}$. $\Box$

\begin{rem} \label{r1}
Let $\mathbf{u}$ be a $TB$-sequence in a countable group $G$. By proposition 3.2 \cite{DMT}, the weight $\omega (G, \sigma_{\mathbf{u}})$ of $(G, \sigma_{\mathbf{u}})$ is $| s_{\mathbf{u}} (G^{\wedge})|$. Since $s_{\mathbf{u}} (G^{\wedge})$ is polishable, $\omega (G, \sigma_{\mathbf{u}})$ is either $\mathfrak{c}$ or $\aleph_0$. $\Box$
\end{rem}

The following lemma plays an important role for the proofs of corollary \ref{c2}.
\begin{lem} \label{l3}
{\it Let $H$ be a dually closed and dually embedded subgroup of a topological group $G$. Then $\mathbf{n} (H) =\mathbf{n} (G)$. }
\end{lem}

\pr Let $x\in \mathbf{n} (H)$. If $\chi \in G^{\wedge}$, then $\chi |_H \in H^{\wedge}$. Hence $(\chi , x) = (\chi |_H , x) =1$ and $x\in \mathbf{n} (G)$.

Let $x\in \mathbf{n} (G)$. Then $x\in H$ (if not, then, since $H$ is dually closed, there exists $\chi \in G^{\wedge}$ such that $(\chi , x)\not= 1$). Let $\eta \in H^{\wedge}$. Then there exists $\chi \in G^{\wedge}$ such that $\chi |_H =\eta $ ($H$ is dually embedded). Then $(\eta, x) = (\chi , x)=1$ and $x\in \mathbf{n} (H)$. $\Box$
\begin{rem} \label{r2}
$\mathbf{n} (G)$ can be characterized as follows: $\mathbf{n} (G)$ is the unique  maximal dually closed subgroup $H$ of $G$ such that there is no non-zero character of $H$ which can be extended to a character of $G$.
\end{rem}

{\bf Proof of corollary \ref{c2}.} The first assertion follows from lemma \ref{l1}. Let us prove the second one. By the definition of $T$-sequence, the subgroup $\langle \mathbf{u}\rangle $ is open in $(H,\mathbf{u})$. Since every open subgroup is dually closed and dually embedded \cite{Nob}, by lemma \ref{l3}, we have $\mathbf{n} (H,\mathbf{u}) = \mathbf{n} (\langle \mathbf{u}\rangle, \mathbf{u})$. If $(H,\mathbf{u})$ is minimally almost-periodic, we have $H=\mathbf{n} (H,\mathbf{u}) = \mathbf{n} (\langle \mathbf{u}\rangle, \mathbf{u}) \subset \langle \mathbf{u}\rangle$. Hence $H=\langle \mathbf{u}\rangle$. Since $\langle \mathbf{u}\rangle$ is countable, the assertion follows from theorem \ref{t3}. $\Box$

{\bf Proof of proposition \ref{p1}.}  Since $(\widehat{X}, \mathbf{u})$ is a $k$-space (corollary 4.1.5 \cite{ZP2}), then $\alpha = \alpha_{(\widehat{X}, \mathbf{u})}$ is continuous (corollary 5.12 \cite{Aus}).

1) Since $(\widehat{X}, \mathbf{u})$ is complete (theorem 2.3.11 \cite{ZP2}), by proposition 6.12 \cite{Aus}, $\alpha$ is an embedding with the closed image. Since $\alpha (X^{\wedge})$ is dense, $\alpha$ is a topological isomorphism.

Since $\alpha$ is bijective, then $\mathbf{n} (\widehat{X}, \mathbf{u}) =\{ 0\}$ and $s_{\mathbf{u}} (X)$ is dense in $X$ by theorem \ref{t3}.

2) Let $i: G_{\mathbf{u}} \to X$ be a continuous monomorphism. Since $i^{\ast} = \alpha$ algebraically, $i^{\ast} (X^{\wedge})$ is dense. Thus $i^{\ast\ast} :  G_{\mathbf{u}}^{\wedge\wedge} \to X$ is a monomorphism. As it was mentioned,  $\alpha : (\widehat{X}, \mathbf{u}) \to G_{\mathbf{u}}^{\wedge}$ is continuous. Hence, $u_n \to 0$ in $G_{\mathbf{u}}^{\wedge}$. If $x\in
G_{\mathbf{u}}^{\wedge\wedge}$, then $(u_n, x) \to 0$. Hence $\alpha_{G_{\mathbf{u}}}$ is bijective.
Since $G_{\mathbf{u}}$ is Polish, then $\alpha_{G_{\mathbf{u}}}$ is a topological isomorphism \cite{Cha}. $\Box$

{\bf Proof of theorem \ref{t6}.} $\Rightarrow$ Let $(\widehat{X}, \mathbf{u})$ be reflexive. Then $\mathbf{n} (\widehat{X}, \mathbf{u})=0$. Hence, by theorem \ref{t3}, ${\rm Cl}_X (s_{\mathbf{u}} (X)) =\mathbf{n} (\widehat{X}, \mathbf{u})^{\perp} =X$ and $s_{\mathbf{u}} (X)$ is dense.

By theorem \ref{t2}, we have
$$
G_{\mathbf{u}}^{\wedge} = (\widehat{X}, \mathbf{u})^{\wedge\wedge} =(\widehat{X}, \mathbf{u}) = \left( X^{\wedge} \times \mathbb{Z}_0^{\infty} \right) / G_{\mathbf{u}}^{\perp}.
$$
Thus any $\chi\in G_{\mathbf{u}}^{\wedge}$ can be extended to some $\widetilde{\chi} \in X^{\wedge} \times \mathbb{Z}_0^{\infty} =\left( X\times \mathbb{T}^H_0 \right)^{\wedge}$. Thus $G_{\mathbf{u}}$ is dually embedded in $X\times \mathbb{T}^H_0$.

$\Leftarrow$ Let $\chi \in G_{\mathbf{u}}^{\wedge}$. Since $G_{\mathbf{u}}$ is dually embedded and $X\times \mathbb{T}^H_0$ is reflexive, then $\alpha = \alpha_{(\widehat{X}, \mathbf{u})} : (\widehat{X}, \mathbf{u})\to (\widehat{X}, \mathbf{u})^{\wedge\wedge} =G_{\mathbf{u}}^{\wedge}$ is a continuous bijection. Hence $(\widehat{X}, \mathbf{u})$ is reflexive by proposition \ref{p1}. $\Box$

{\large Department of Mathematics, Ben-Gurion University of the
Negev,

Beer-Sheva, P.O. 653, Israel}

{\it E-mail address}: $\quad$ saak@math.bgu.ac.il

\end{document}